\documentclass[a4paper,12pt]{amsart} 
\usepackage[cp866]{inputenc}
\usepackage[T2A]{fontenc}
\usepackage[russian]{babel}
\usepackage{amsfonts, amssymb, amsmath, wrapfig}
\usepackage{longtable}
\usepackage[all]{xy}
\usepackage{epsfig}

\theoremstyle{plain}
\newtoks\thehProclaim
\newtheorem*{Proclaim}{\the\thehProclaim}
\newenvironment{proclaim}[1]{\thehProclaim{#1}\begin{Proclaim}}{\end{Proclaim}}

\def\diag{\text{\rm diag}}

\def\rk{\operatorname{rk}}

\def\GF#1{{\mathbb F}_{\!#1}}

\def\e{\varepsilon}
\def\a{\alpha}
\def\b{\beta}
\def\g{\gamma}

\def\A{\operatorname{A}}

\def\C{\operatorname{C}}

\def\F{\operatorname{F}}

\def\K{\operatorname{K}}

\def\Sp{\operatorname{Sp}}

\def\SL{\operatorname{SL}}
\def\GL{\operatorname{GL}}

\def\St{\operatorname{St}}

\def\SN{\operatorname{SN}}
\def\sr{\operatorname{sr}}
\def\asr{\operatorname{asr}}
\def\tg{\operatorname{tg}}

\def\Int{{\mathbb Z}}
\def\Rat{{\mathbb Q}}
\def\Co{{\mathbb C}}

\begin{document}

\title[Unitriangular factorisations of Chevalley groups]
{unitriangular factorisations of Chevalley groups}
\author{N.~A.~Vavilov, A.~V.~Smolensky, B.~Sury}

\address{St. Petersburg\newline
State University\newline
University pr. 28,
Peterhof, \newline 198504 St. Petersburg, Russia
\newline
Indian Statistics Institute
Bangalore}
\email {nikolai-vavilov\@yandex.ru
\newline
andrei.smolensky\@gmail.com
\newline
surybang\@gmail.com}
 \date{ 27 мая 2011 г.}

\thanks
{The research of the first author was carried out in the
framework of the RFFI projects 09-01-00784, ``Efficient
generation in groups of Lie type'' (PDMI RAS),
09-01-00878 ``Overgroups of reductive groups in algebraic
groups over rings'' (SPbGU), 10-01-90016
``The study of structure of forms of reductive groups and behaviour
of small unipotent elements in representations of algebraic groups''
(SPbGU). The work of the first and the third author was conducted
in the framework of a joint Russian--Indian project RFFI 10-01-92651
``Higher composition laws, algebraic groups,
K-theory and exceptional groups'' (SPbGU). Apart from that, the first
author was partially supported by the RFFI projects 09-01-00762
(Siberian Federal University), 09-01-91333 (POMI RAS) and
11-01-00756 (RGPU). The work of the first and the second author
was supported also by the State Financed task project 6.38.74.2011
at the Saint Petersburg State University
``Structure theory and geometry of algebraic groups, and their
applications in representation theory and algebraic $\K$-theory''.
The third author is grateful to the Saint Petersburg Department of
Steklov Mathematical Institute and to the Saint Petersburg State
University for the invitation to visit Saint Petersburg in
May--June 2011.}

\maketitle

\par
In the present paper we show that comparing results of Hyman Bass
\cite{bass64} with those of Oleg Tavgen \cite{tavgen90b}, one immediately 
gets the following result, which is both more general, and more 
precise than all recent results pertaining to unitriangular 
factorisations.
\begin{proclaim}
{Theorem 1} Let\/ $\Phi$ be a reduced irreducible root system
and\/ $R$ be a commutative ring such that $\sr(R)=1$. Then
the elementary Chevalley group\/ $G(\Phi,R)$ admits
unitriangular factorisation
$$ E(\Phi,R)=U(\Phi,R)U^-(\Phi,R)U(\Phi,R)U^-(\Phi,R). $$
\noindent
of length\/ $4$.
\end{proclaim}

On the other hand, since $U^-(\Phi,R)\cap B(\Phi,R)=1$, one has
$$ T(\Phi,R)\cap U(\Phi,R)U^-(\Phi,R)U(\Phi,R)=1. $$
\noindent
In other words, 1 is the {\it only\/} element of the torus
$T(\Phi,R)$, that admits a unitriagular factorisation of
length $<4$. Thus, if the ring $R$ has at least one non-trivial
unit, factorisation obtained in Theorem 1 is the best possible.
\par
What is truly amazing here, is that the usual {\it linear\/}
stable rank condition works for groups of all types!
Under somewhat stronger stability conditions, a similar result
holds also for twisted Chevalley groups. However, the analysis
of twisted groups requires much more detailed calculations for
groups of Lie ranks 1 and 2, where one cannot simply invoke
known results. The proofs for this case are relegated to the
next paper by the authors.
\par
Clearly, for rings of dimension $\ge 1$ in general there is no
way to obtain unitriangular factorisations of length 4.
However, for some particularly nice rings of dimension 1 one
can obtain unitriangular factorisations of length 5 or 6.
The second main result of the present paper gives a simplest
such example.

\begin{proclaim}{Theorem 2}
Let\/ $\Phi$ be a reduced irreducible root system and\/ $p\in\Int$
be a rational prime. Under assumption of the Generalised Riemann's
Hypothesis the simply connected Chevalley group
$G\Big(\Phi,\Int\Big[{1\over p}\Big]\Big)$
admits unitriangular factorisation
$$ G\bigg(\Phi,\Int\bigg[{1\over p}\bigg]\bigg)=
\bigg(U\bigg(\Phi,\Int\bigg[{1\over p}\bigg]\bigg)
U^-\bigg(\Phi,\Int\bigg[{1\over p}\bigg]\bigg)\bigg)^3 $$
\noindent
of length\/ $6$.
\end{proclaim}

Presently, we are working towards generalisation of this result
to other Hasse domains and plan to return to this topic in a
separate paper.
\par
The first part of the present paper is the term paper of the second
author, under the supervision of the first author, whereas the
second part is a by-product of our joint work on arithmetic problems
our cooperative Russian--Indian project ``Higher composition laws,
algebraic groups, K-theory and exceptional groups'' at the Saint
Petersburg State University, Tata Institute of Fundamental Research
(Mumbai) and Indian Statistical Institute (Bangalore).
\par
Since there are a considerable number of papers addressing
some aspects of unitriangular factorisations, and many of
these papers give no reference to previous works, where
similar --- or even stronger! --- results are established,
we start our paper with a brief survey of the known results,
as we understand them. This is done in \S\S1--4. After that
in \S\S5 and 6 we introduce basic notation. In \S7 we
prove a version of Tavgen's theorem on rank reduction, which
immediately implies Theorem 1. In \S8 we discuss connection of
unitriangular factorisations with Bruhat and Gauss decompositions.
In \S9 we consider Chevalley groups over arithmetic rings and
prove Theorem 2. Finally, in \S10 we mention some further related
topics and state several unsolved problems.


\section*{\S1. Existence of unitriangular factorisations}

Recently, the following problem cropped up anew, in several
independent contexts: find the shortest factorisation
$G=UU^-UU^-\ldots U^{\pm}$ of a Chevalley group $G=(\Phi,R)$,
in terms of the unipotent radical $U=U(\Phi,R)$ of the
standard Borel subgroup $B=B(\Phi,R)$, and the unipotent
radical $U^-=U^-(\Phi,R)$ of the opposite Borel subgroup
$B^-=B^-(\Phi,R)$.
\par
Let us mention three such extensive clusters of various subjects,
where it occurred. In the following sections we see that although
information between various subjects within one cluster is
transmitted with serious delays, still communication is {\it
incredibly\/} faster and more effective, than information transfer
between clusters themselves.

\smallskip
$\bullet$ Algebraic $K$-theory, structure theory of algebraic
groups, theory of arithmetic groups.

\smallskip
$\bullet$ Theory of finite and profinite groups, asymptotic
group theory, finite geometries.

\smallskip
$\bullet$ Computational linear algebra, wavelet theory,
computer graphics, control theory.

\smallskip
Let us list some existing results concerning factorisations
of the form $G=UU^-UU^-\ldots U^{\pm}$. Firstly, one has to
establish {\it existence\/} of such factorisations, and,
secondly, estimate their {\it length\/}.
\par
Clearly, existence of some unitriangular factorisation is
equivalent to two following conditions.

\smallskip
$\bullet$ The Chevalley group $G(\Phi,R)$ coincides with its
elementary subgroup $E(\Phi,R)$, spanned by elementary
generators.

\smallskip
$\bullet$ The width of the elementary Chevalley group $E(\Phi,R)$
with respect to elementary generators is bounded.

\smallskip
The answer to both of these questions is, in general, negative
beyond retrieve, so that one can only expect existence of
unitriangular factorisations for some very special classes of
rings. There is {\it vast\/} literature, dedicated to both of
these problems, and we will not even to attempt to address
them in full. Instead, we confine ourselves to a brief excerpt
of the introduction to \cite{SV}, referring to \cite{tavgen92,
morris} for a broader picture and further references.

\smallskip
$\bullet$ In the case when the ring $R$ is semilocal ---
for instance, when it is of finite dimension over a field ---
for simple connected Chevalley groups one has the equality
$E(\Phi,R)=G(\Phi,R)$. Bounded width in terms of elementary
generators immediately follows from Gauss decomposition.

\smallskip
$\bullet$ It is classically known that expressions of a matrix in
the group $\SL(2,R)$ over a Euclidean ring $R$ as a product of
elementaries correspond to continued fractions. Existence of
arbitrary long division chains in $\Int$ shows that the group
$\SL(2,\Int)$ cannot have bounded width in elementary generators.
Further such results are discussed in \cite{cohn,CW}.

\smallskip
$\bullet$ On the other hand, David Carter and Gordon Keller 
\cite{CK83,CK84} and Oleg Tavgen \cite{tavgen90a,tavgen90b,tavgen92} 
demonstrated that for a Dedekind ring of
arithmetic type $R$, the simply connected Chevalley groups
$G(\Phi,R)=E(\Phi,R)$ of rank $\ge 2$ have bounded width in
elementary generators.

\smallskip
$\bullet$ Wilberd van der Kallen \cite{vdK82} made a striking discovery 
that in general even Chevalley groups of rank $\ge 2$ over a Euclidean
ring may have infinite width in elementary generators. More
precisely, he proved that $\SL(3,\Co[t])$ -- and, since
$\sr(\Co[t])=2$, all groups $\SL(n,\Co[t])$ for $n\ge 4$ -- do not
have bounded width in elementary generators.

\smallskip
Thus, in general one can hope to establish existence of
unitriangular factorisations only over some very special
rings of dimension $\le 1$.
\par
In the next three sections we list some works addressing the length
of unitriangular factorisations. Already this short summary shows
that experts in one field are usually completely unaware of the
standard notions and results in another field. It is hard to
imagine, how much time and energy could have been saved, should
millions of programmers and engineers learn the words {\it parabolic
subgroup\/} and {\it Levi decomposition\/}, rather than persisting
in retarded matrix manipulations.
\par
The same remark applies to Application~1.1 of the work \cite{BNP}, 
where, at least for the linear case, the result was known for quite 
some time, in a much larger generality, and with better bound, to 
experts in other subjects, both in algebraic $K$-theory \cite{DV88} 
and in computational linear algebra \cite{strang}. Obviously, after 
its driving forces are understood for the case of $\SL(n,R)$, its 
generalisation to Chevalley groups does not require any serious 
intellectual effort, just the mastery of the appropriate techniques.

\section*{\S2. The length of factorisations: linear groups}

For the linear case a systematic study of this problem was
started in a paper by Keith Dennis and Leonid Vaserstein \cite{DV88}.
Originally, their interest to this problem was motivated
by the following observation.

\smallskip
$\bullet$ Any element of the group $U(n,R)$, $n\ge 3$, over
an arbitrary associative ring $R$ is a product of not more
than two commutators in the elementary group $E(n,R)$,
Lemma 13. Earlier, van der Kallen cite{vdK82} noticed that the
elements of $U(n,R)$ are products of not more than three
commutators.

\smallskip
Let us state several typical results of \cite{DV88}, pertaining
to the length of unitriangular factorisations.

\smallskip
$\bullet$ If for a ring $R$ the group $E(m,R)$, $m\ge 2$,
is represented as a finite product $UU^-UU^-\ldots U^{\pm}$
with $L$ factors, then all groups $E(n,R)$, $n>m$,
can be presented as the product of the same form, with
the same number of factors, Lemma~7.

\smallskip
$\bullet$ If the stable rank of the ring $R$ equals 1, then,
as is classically known from the work of Hyman Bass \cite{bass64},
the group $E(2,R)$ -- and thus by the previous item all
groups $E(n,R)$ -- admit unitriangular factorisation
$$  E(n,R)=U(n,R)U^-(n,R)U(n,R)U^-(n,R) $$
\noindent
of length 4.

\smallskip
To the authors of \cite{DV88}, as to all experts in algebraic
$K$-theory, this fact is obvious to such an extent, that
in \cite{DV88} it is not even separately stated, and just casually
mentioned inside the proof of Theorem~6. Nevertheless, since
this is one of the key steps in the proof of Theorem~1, at
that the {\it only one\/}, which invokes stability condition,
we reproduce its proof, which is an adaptation of a more
general similar calculation, implemented in \cite{bass64}.
\par
Recall that a ring $R$ has {\bf stable rank $1$\/}, if for
all $x,y\in R$, which generate $R$ as a {\it right\/} ideal,
there exists a $z\in R$ such that $x+yz$ is {\it right\/}
invertible. In this case we write~$\sr(R)=1$.
\par
The simplest and the most characteristic example of a ring
of stable rank 1 are semilocal rings. However, there are many
far less trivial examples, such as the ring of all algebraic
integers. Many further examples and further references can be
found in \cite{vaser}.
\par
It is classically known that rings of stable rank 1 are actually
{\it weakly finite\/} (Kaplansky---Lenstra theorem), so that in
their definition one could from the very start require
that $x+yz\in R^*$. Since for the linear case the result is
well known, and Chevalley groups of other types only exist over
commutative rings, from here on we assume that the ring $R$ is
commutative, in which case the proof below at the same time
demonstrates that $\SL(2,R)=E(2,R)$.
\par
In the following proof, and henceforth, in discussing linear case
we use standard matrix notation. In particular, $e$ denotes
identity matrix, $e_{ij}$, $1\le i,j\le n$, denotes a standard
matrix unit, i.~e.~the matrix which has 1 in the position $(i,j)$
and zeroes elsewhere. Further, for $1\le i\neq j\le n$ and
$\xi\in R$ one denotes by $t_{ij}(\xi)$ the elementary
transvection $e+\xi e_{ij}$. 
The matrix entry of $g$ in position $(i,j)$ will be denoted by
$g_{ij}$, whereas the corresponding entry of the inverse matrix
$g^{-1}$ will be denoted by $g_{ij}'$.

\begin{proclaim}
{Lemma 1} Let\/ $R$ be a commutative ring of stable rank\/ $1$. Then
$$ \SL(2,R)=U(2,R)U^-(2,R)U(2,R)U^-(2,R). $$
\end{proclaim}

\begin{proof}[Proof] Let us trace, how many elementary transformations
one needs to bring an arbitrary matrix
$g=\begin{pmatrix} a&b\\ c&d\\ \end{pmatrix}\in\SL(2,R)$
to the form $e$. We will not introduce new notation at each step,
but rather replace the matrix $g$ by its current value, as is
common in computer science. Obviously, its entries $a,b,c,d$
should be also reset to their current values at each step.

\smallskip\noindent
{\bf Step 1.} \ Multiplication by a single {\it lower\/}
elementary matrix on the right allows to make the element in
the South-West corner invertible. Indeed, since the rows of
the matrix are unimodular, one has $cR+dR=R$. Since $\sr(R)=1$,
there exists such an $z\in R$, that $c+dz\in R^*$. Thus,
$$ gt_{21}(z)=\begin{pmatrix} a+bz&b\\ c+dz&d\\
\end{pmatrix}, $$
\noindent
where $c+dz\in R^*$.

\smallskip\noindent
{\bf Step 2.} \ Thus, we can assume that $c\in R^*$.
Multiplication by a single {\it upper\/} elementary matrix
on the right allows to make the element in the South-East
corner equal to 1. Indeed,
$$ gt_{12}(c^{-1}(1-d))=\begin{pmatrix} a&b+ac^{-1}(1-d)\\
c&1\\
\end{pmatrix}, $$

\smallskip\noindent
{\bf Step 3.} \ Thus, we can now assume that $d=1$. Multiplication
by a single {\it lower\/} elementary matrix on the right allows to
make the element in the South-West corner equal to 0. Indeed,
$$ gt_{21}(-c)=\begin{pmatrix} a-bc&b\\ 0&1\\
\end{pmatrix}. $$
\par
Since $\det(g)=1$, the matrix on the right hand side is
equal to $t_{12}(b)$. Bringing all elementary factors to the right
hand side, we see that any matrix $g$ with determinant 1 can be
expressed as a product of the form $t_{12}(*)t_{21}(*)t_{12}(*)t_{21}(*)$,
as claimed.
\end{proof}

Now, let us return to the Dennis---Vaserstein paper \cite{DV88}.

\smallskip
$\bullet$ If $R$ is a Boolean ring, in other words, if $x^2=x$
for all $x\in R$, then the elementary group $E(n,R)$, $n\ge 2$,
admits the unitreiangular factorisation
$$ E(n,R)=U(n,R)U^-(n,R)U(n,R) $$
\noindent
of length 3. For commutative rings the converse is obviously true:
if for some $n\ge 2$ the group $E(n,R)$ admits unitriangular
factorisation of length 3, then the ring $R$ is Boolean, Remark~15.

\smallskip
$\bullet$ Let $d=\diag(\e_1,\ldots,\e_n)$, $\e_i\in R^*$,
be a diagonal matrix such that $\e_1\ldots\e_n=1$.
Then $d\in U(n,R)U^-(n,R)U(n,R)U^-(n,R)$, Lemma 18.
The proof in \cite{DV88} is based on general position arguments.
In \S8 we give another proof of a slightly stronger
fact, based on induction on rank, which can be easily
generalised to all Chevalley groups.

\smallskip
$\bullet$ If the stable rank of the ring $R$ is finite for
some $m\ge 2$, the group $E(m,R)$ has finite width with
respect to the elementary generators, then for all sufficiently
large $n$ one has the unitriangular factorisation
$$ E(n,R)=(U(n,R)U^-(n,R))^3 $$
\noindent
of length 6.

\smallskip
$\bullet$ It is extremely suggestive to compare this result with
{\it Sharpe decomposition}. Recall that it is established in 
\cite{sharpe} that for any associative ring $R$ the finitary 
elementary group $E(R)=\lim\limits_{\longrightarrow}E(n,R)$, 
$g\mapsto g\oplus 1$, admits decomposition
$$
\begin{aligned}
E(R)&=B(R)N(R)U(R)U^-(R)\\
&=B(R)U^-(R)U(R)U^-(R)=
U(R)D(R)U^-(R)U(R)U^-(R),
\end{aligned}
$$
\noindent where $B(R)$, $N(R)$, $U(R)$, $U^-(R)$, $D(R)$ denote
inductive limits of the corresponding groups of finite degree, with
respect to the same embedding. Since $D(R)\subseteq
U(R)U^-(R)U(R)U^-(R)$, it follows that for an arbitrary associative
ring the finitary elementary group $E(R)$ admits unitriangular
factorisation
$$ E(R)=(U(R)U^-(R))^3 $$
\noindent
of length 6.

\smallskip
$\bullet$ Recently, motivated by applications in factorisation
of integer matrices, Thomas Laffey cite{laffey-a,laffey-b} explicitly
calculated the bound in Den\-nis---Vaserstein theorem for the
ring $R=\Int$. In particular, unitriangular factorisation
$$ \SL(n,\Int)=(U(n,\Int)U^-(n,\Int))^3 $$
\noindent
of length 6 holds for any $n\ge 82$.

\section*{\S3. The length of factorisations: Chevalley groups}

The proof that Chevalley groups over Dedekind rings of arithmetic
type have finite width in elementary generators, given by Oleg
Tavgen \cite{tavgen90b}, relied on reduction to groups of rank 2. 
The base of this reduction is the following fact, see 
\cite{tavgen90b}, Proposition~1.  

\smallskip
$\bullet$ If all Chevalley groups of a certain rank $l$
admit unitriangular factorisation $G=(UU^-)^L$ of certain
length, then all groups of larger ranks also admit unitriangular
factorisations of the same length. Actually, our Theorem 3
is a minor elaboration of Tavgen's idea.

\smallskip
$\bullet$ Another elaboration of the same idea can be found
in a recent paper by Andrei Rapinchuk and Igor Rapinchuk \cite{RR}.
There, the same idea is used to establish that Chevalley
groups over a local ring $R$ can be factored as
$$ G(\Phi,R)=\big(U(\Phi,R)U^-(\Phi,R)\big)^4. $$

\par
$\bullet$ For symplectic group our Theorem~1 was established 
by You Hong \cite{you}. For orthogonal groups a similar, but
weaker result, under stronger stability conditions was established
by Frank Arlinghaus, Leonid VAserstein and You Hong \cite{AVY}.

\par\smallskip 
Over fields, especially over finite fields, there was a
vivid interest in explicit calculation of length, over the
last decade. Let $\GF{q}$, $q=p^m$, be a finite field of
characteristic $p$. In this case the group
$U(\Phi,q)=U(\Phi,\GF{q})$ is a Sylow $p$-subgroup of the Chevalley
group $G(\Phi,q)=G(\Phi,\GF{q})$. Thus, in this case calculation
of the minimal length of unitriangular factorisation is
essentially equivalent to the calculation of the minimal length
of factorisation of a finite simple group of Lie type
in terms of its Sylow $p$-subgroups, in the defining characteristic.

\smallskip
$\bullet$ Martin Liebeck and Laszlo Pyber \cite{LP}, Theorem D,
prove that finite Chevalley groups admit unitriangular factorisation
$$ G(\Phi,q)=(U(\Phi,q)U^-(\Phi,q))^{6}U(\Phi,q) $$
\noindent
of length 13. Actually, they also consider twisted Chevalley groups
and obtain for them a similar bound, with the only exception
of the senior Ree groups ${}^2\!\F_4(q)$, for which they prove
existence of factorisation $G=(UU^-)^{12}U$.

\smallskip
$\bullet$ Laszlo Babai, Nikolay Nikolov and Laszlo Pyber \cite{BNP},
Application~1.1, prove that finite Chevalley groups admit
unitriangular factorisation
$$ G(\Phi,q)=U(\Phi,q)U^-(\Phi,q)U(\Phi,q)U^-(\Phi,q)U(\Phi,q) $$
\noindent of length 5. They also obtain a similar result with the
same bound for twisted groups.
\par
As opposed to \cite{LP}, the proofs in the paper \cite{BNP} are of 
extremely mysterious nature. In the final count, they rely on the 
product growth estimates of the following sort. Let 
$X,Y\subseteq G$ be two nonempty subsets of a finite group $G$. Then
$$ |XY|\ge\min\left({|G|\over 2},{m|X|\cdot|Y|\over 2|G|}\right), $$
\noindent where $m$ denotes the smallest dimension of a nontrivial
real representation of the group $G$. In particular, in the case
when the group $G$ is simple, and the orders $|X|,|Y|$ are large
enough, but still much smaller than $|G|$, products rapidly grow:
the order $|XY|$ is much larger, than $\max(|X|,|Y|)$.
\par
In this connection, the authors of \cite{BNP} make the following 
extravagant claim:
``For the most part we can argue by the size of certain subsets,
ignoring the structure and thus greatly simplifying the proofs and
at the same time obtaining considerably better, nearly optimal
bounds.''
\par
Frankly, the proposal to completely ignore the structure of the
objects in question, seems to us excessive. Asymptotic methods are
in their place in the study of profinite groups, or groups of
infinite rank. On the other hand, when applied to the finite groups
proper, or to groups of Lie type of finite rank, asymptotic methods
should be considered as {\it surrogate\/} of structural algebraic
methods, as a method to verify the result we are interested in,
which is {\it identical\/} to a real proof.
\par
Certainly, in many situations such a genuine algebraic proof does
not exist. What is worse, for many interesting problems such
algebraic solutions are not in sight, since they would either
require {\it enormous\/} case by case analysis, or cannot be
obtained by the present day methods.
\par
As an example of the first situation, we can cite the remarkable
paper by Martin Liebeck and Aner Shalev \cite{LS96}, where it is
proven that --- apart from three exceptional series, Suzuki
groups and groups $\Sp(4,2^m)$ and $\Sp(4,3^m)$ --- almost all
finite simple groups are $(2,3)$-generated. For each individual
group of Lie type, and at that not only over a finite field,
but also over finitely generated rings, it is basically clear
how to prove its $(2,3)$-generation, or the lack thereof.
In fact, for important classes of groups, including all groups
of sufficiently large Lie ranks (in the order of magnitude of
several dozens, for finite simple groups) such constructive
algebraic proofs have been obtained. However, a complete
analysis of groups of small ranks requires computational effort
so extensive that it has not been completed yet, despite
enormous work by many authors, over more than two decades.
\par
As an example of the second situation, we could cite recent
papers on the verbal width of finite simple groups. An amazing
general result by Aner Shalev cite{shalev09}, asserts that for {\it any\/}
nontrivial word $w\neq 1$ the verbal width of almost all
finite simple groups $G$ equals 3. In other words, every element
of the group $G$ can be presented as a product of not more
than three values of the word $w$. See also \cite{GM,LS,LS01,LP,shalev07},
where one can find further results in the same style, improvement
of the estimate for verbal width to 2, in some cases, and
references to preceding works. As far as we know, similar
results are not available not to say over rings, but even over
infinite fields. What is worse, such similar results are not
available not just for all words, but for simplest {\it specific\/}
words, such as powers. What is still worse, there is no clear
understanding, how one could prove such results at all.
\par
Nevertheless, we believe that whenever one can apply algebraic
methods, coming from the structure theory and representation
theory of algebraic groups, they would {\it invariably\/}
give stronger and more general results. For those results that
hold over arbitrary fields they should also give better bounds,
and {\it hopefully\/}, have simpler proofs.

\section*{\S4. The length of factorisations:\\
computational linear algebra}

Definitely, unitriangular factorisations for the group $\SL(n,K)$
are so natural and obvious, that there is little doubt that they
should have been known to experts in linear algebra for quite some
time. This is however the earliest reference we could trace.

\smallskip
$\bullet$ Gilbert Strang \cite{strang} noticed that all groups 
$\SL(n,K)$ over a field $K$ admit unitriangular factorisation
$$ \SL(n,K)=U^-(n,K)U(n,K)U^-(n,K)U(n,K) $$
\noindent of length 4. This is what experts in computational linear
algebra call LULU-factorisation, where mnemonically L should be
interpreted as the first letter of the word `lower', whereas U
should be interpreted as the first letter of the word `upper'.

\smallskip
Observe, that this fact got to the household of linear algebra
decades after it became standard in algebraic $K$-theory, and
at that in a much larger generality. However, as we know, the
walls between different branches of mathematics are high.
\par
As another amusing circumstance, we could mention that, as far as we
know, this result was only observed in linear algebra due to
applications in computer graphics! Let us explain in somewhat more
detail, how it happened. It is hard for us to construe the {\it
precise\/} sense in which experts in computer algebra use the term
{\it shear\/}. As a first approximation one can think of shear as an
arbitrary unipotent element of $\SL(n,K)$. Or, at least, majority of
authors in this field call arbitrary elements of the groups $U(n,K)$
and $U^-(n,K)$ shears.
\par
However, there is no doubt how one should interpret the term
{\it one-dimensional shears\/}. These are precisely transvections,
not necessarily elementary. An exclusive role is played by
transvections concentrated in one row --- they are called
{\it beam shears\/} --- or in one column --- they are called
{\it slice shears\/}. The point is that at the level of pixels
such transvections correspond to string copy with offset, which
admits {\it extremely\/} efficient hardware implementations.
\par

\def\phi{\varphi}
Alain Paeth cite{paeth} proposed the following cute method to 
implement a 2D rotation:
$$ \begin{pmatrix}
\cos(\phi)&\sin(\phi)\\
-\sin(\phi)&\cos(\phi)\\
 \end{pmatrix}=
 \begin{pmatrix}
1&\tg(\phi/2)\\
0&1\\
 \end{pmatrix}
 \begin{pmatrix}
1&0\\
-\sin(\phi)&1\\
 \end{pmatrix}
 \begin{pmatrix}
1&\tg(\phi/2)\\
0&1\\
\end{pmatrix}. $$
\noindent
As we already mentioned, the shears themselves are implemented
at the level of data transfer, this method turned out to be much
more efficient than coordinate conversion.
\par
At this point it is natural to ask, whether 3D rotations admit a
similar efficient hardware implementation? Initially, one attempted
to first decompose a 3D rotation into a product of three 2D
rotations, and then each of them into a product of three
transvections. However, shortly thereafter Tommaso Toffoli and Jason
Quick \cite{TQ} proposed a scheme based upon decomposition of a 3D
rotation into a product of three unipotent matrices, and designed
corresponding hardware.
\par
Let us recall, that a three-dimensional rotation $g$ is completely
determined by its Euler angles $(\a,\b,\g)$, for instance, as follows:
{\small 
$$ g= \begin{pmatrix}
c(\a)c(\b)c(\g)-s(\a)s(\g)
&-c(\a)c(\b)s(\g)-s(\a)c(\g)
&c(\a)s(\b)\\
\noalign{\vskip 5truept}
s(\a)c(\b)c(\g)+c(\a)s(\g)
&-s(\a)c(\b)s(\g)+c(\a)c(\g)
&s(\a)s(\b)\\
\noalign{\vskip 5truept}
-s(\b)c(\g)&s(\b)s(\g)&c(\b)\\
 \end{pmatrix},
 $$
 }
 \\
\noindent
where $c(\varphi)$ and $s(\varphi)$ denote $\cos(\varphi)$
and $\sin(\varphi)$, respectively. It is easy to notice
--- this was the starting point of \cite{TQ} --- that $g$
admits the following unitriangular factorisation of length 3:
 \begin{multline*}
g=
 \begin{pmatrix}
1&-\displaystyle{\tg\left({\a+\g\over 2}\right)}&
\displaystyle{\cos(\a)\tg\left({\b\over 2}\right)}\\
\noalign{\vskip 5truept}
0&1&\displaystyle{\sin(\a)\tg\left({\b\over 2}\right)}\\
\noalign{\vskip 5truept}
0&0&1\\
 \end{pmatrix}   \\
 \times\begin{pmatrix}
1&0&0\\
\noalign{\vskip 5truept}
\sin(\a+\g)&1&0\\
\noalign{\vskip 5truept}
-\cos(\g)\sin(\b)&-
\displaystyle{\displaystyle{\sin\left({\a-\g\over 2}\right)}
\over\displaystyle{\cos\left({\a+\g\over 2}\right)}}
\sin(\b)&1\\
 \end{pmatrix}  \\
 \times\begin{pmatrix}
1&-\displaystyle{\tg\left({\a+\g\over 2}\right)}&
\displaystyle{\displaystyle{\cos\left({\a-\g\over 2}\right)}\over
\displaystyle{\cos\left({\a+\g\over 2}\right)}}\tg\left({\b\over 2}\right)\\
\noalign{\vskip 5truept}
0&1&-\displaystyle{\sin(\g)\tg\left({\b\over 2}\right)}\\
\noalign{\vskip 5truept}
0&0&1\\
 \end{pmatrix}.
\end{multline*}

\par
The work by Strang \cite{strang}, cited at the beginning of the section,
emerged as an attempt to generalise this formula to the case of an
arbitrary $n$. After that many subsequent papers appeared, which
discussed various aspects of such decompositions. To convey some
flavour of this activity, let us state a couple of typical results
from the next paper by Toffoli \cite{toffoli}.

\smallskip
$\bullet$ Almost all elements $\SL(n,{\mathbb R})$ admit unitriangular
factorisations of length 3 --- ULU-factorisation, as it is called
by the experts in computational linear algebra. Obviously, as
opposed to the previous section, here the expression `almost all'
should be interpreted in terms of Lebesgue measure.

\smallskip
$\bullet$ For any matrix $g\in\SL(n,K)$ there exists such a permutation
matrix $(\pi)$, $\pi\in\operatorname{S}_n$, that at least one of the
matrices $(\pi)g$ or $(\pi)^{-1}g$ admits unitriangular factorisation
of length 3

\smallskip
We will make no attempts to systematically cover subsequent literature
in the field, and limit ourselves to several typical works
\cite{CK,hao,LHW,SH}, where one can find further references.

\section*{\S5. Some subgroups of Chevalley groups}

\def\wP{\mathcal P}
\def\wQ{\mathcal Q}

Our notation pertaining to Chevalley groups are utterly standard
and coincide with the ones used in \cite{vavilov,VP},
where one can find many further references.
\par
Let $\Phi$ be a reduced irreducible root system of rank $l$,
$W=W(\Phi)$ be its Weyl group and $\wP$ be a weight lattice
intermediate between the root lattice $\wQ(\Phi)$ and the weight
lattice $\wP(\Phi)$. Further, we fix an order on $\Phi$ and denote
by $\Pi=\{\a_1,\ldots,\a_l\}$, $\Phi^+$ and $\Phi^-$ the corresponding
sets of fundamental, positive and negative roots, respectively.
Our numbering of the fundamental roots follows Bourbaki.
Finally, let $R$ be a commutative ring with 1, as usual,
$R^*$ denotes its multiplicative group.
\par
It is classically known that with these data one can associate
the Chevalley group $G=G_{\wP}(\Phi,R)$, which is the group of
$R$-points of an affine groups scheme $G_{\wP}(\Phi,-)$, known as
the Chevalley---Demazure group scheme. Since mostly our results
do not depend on the choice of the lattice $\wP$, in the sequel we
usually assume that $\wP=\wP(\Phi)$ and omit any reference to
$\wP$ in the notation. Thus, $G(\Phi,R)$ will denote the
simply connected Chevalley group of type $\Phi$ over $R$.
\par
In what follows, we fix a split maximal torus $T(\Phi,-)$
of the group scheme $G(\Phi,-)$ and set $T=T(\Phi,R)$.
As usual, $X_{\a}$, $\a\in\Phi$, denotes a unipotent root subgroup
in $G$, elementary with respect to $T$. We fix isomorphisms
$x_{\a}:R\mapsto X_{\a}$, so that
$X_{\a}=\{x_{\a}(\xi)\mid\xi\in R\}$, which are interrelated
by the Chevalley commutator formula, see \cite{steinberg,carter,VP}. 
Further,
$E(\Phi,R)$ denotes the elementary subgroup of $G(\Phi,R)$,
generated by all root subgroups $X_{\a}$, $\a\in\Phi$.
\par
In the sequel the elements $x_{\a}(\xi)$ are called root unipotents
Now, let $\a\in\Phi$ and $\e\in R^{*}$. As usual, we set
$h_{\a}(\e)=w_{\a}(\e)w_{\a}(1)^{-1}$, where
$w_{\a}(\e)=x_{\a}(\e)x_{-\a}(-\e ^{-1})x_{\a}(\e)$. The elements
$h_{\a}(\e)$ are called semisimple root elements. For a simply
connected group one has
$$ T=T(\Phi,R)=\langle h_{\a}(\e),\ \a\in\Phi,\ \e\in R^*\rangle. $$
\noindent
Finally, let $N=N(\Phi,R)$ be the algebraic normaliser of the
torus $T=T(\Phi,R)$, i.~e.\ the subgroup, generated by $T=T(\Phi,R)$
and all elements $w_{\a}(1)$, $\a\in\Phi$. The factor-group
$N/T$ is canonically isomorphic to the Weyl group $W$, and for each
$w\in W$ we fix its preimage $n_{w}$ in $N$.
\par
The following result is obvious, well known, and very useful.

\begin{proclaim}{Lemma 2}
The elementary Chevalley group\/ $E(\Phi,R)$ is generated by
unipotent root elements\/ $x_{\a}(\xi)$, $\a\in\pm\Pi$, $\xi\in R$,
corresponding to the fundamental and negative fundamental roots.
\end{proclaim}

\begin{proof}[Proof]
Indeed, every root is conjugate to a fundamental root by an
element of the Weyl group, while the Weyl group itself is
generated by the fundamental reflections $w_{\a}$,
$\a\in\Pi$. Thus, the elementary group $E(\Phi,R)$ is generated
by the root unipotents $x_{\a}(\xi)$, $\a\in\Pi$,
$\xi\in R$, and the elements $w_{\a}(1)$, $\a\in\Pi$. It remains
only to observe that $w_{\a}(1)=x_{\a}(1)x_{-\a}(-1)x_{\a}(1)$.
\end{proof}

Further, let $B=B(\Phi,R)$ and $B^-=B^-(\Phi,R)$ be a pair of
opposite Borel subgroups containing $T=T(\Phi,R)$, standard
with respect to the given order. Recall that $B$ and $B^-$
are semidirect products $B=T\rightthreetimes U$ and
$B^-=T\rightthreetimes U^-$, of the torus $T$ and their unipotent
radicals
$$
\begin{aligned}
U&=U(\Phi,R)=
\big\langle x_\a(\xi),\ \a\in\Phi^+,\ \xi\in R\big\rangle, \\
\noalign{\vskip 4pt}
U^-&=U^-(\Phi,R)=
\big\langle x_\a(\xi),\ \a\in\Phi^-,\ \xi\in R\big\rangle. \\
\end{aligned}
$$
\noindent
Here, as usual, for a subset $X$ of a group $G$ one denotes by
$\langle X\rangle$ the subgroup in $G$ generated by $X$.
Semidirect product decomposition of $B$ amounts to saying that
$B=TU=UT$, and at that $U\trianglelefteq B$ and $T\cap U=1$.
Similar facts hold with $B$ and $U$ replaced by $B^-$ and $U^-$.
Sometimes, to speak of both subgroups $U$ and $U^-$ simultaneously,
we denote $U=U(\Phi,R)$ by $U^+=U^+(\Phi,R)$.
\par
Generally speaking, one can associate a subgroup $E(S)=E(S,R)$
to any closed subset $S$ in $\Phi$. Recall that a subset $S$
in $\Phi$ is called {\it closed\/}, if for any two roots
$\a,\b\in S$ the fact that $\a+\b\in\Phi$, implies that already
$\a+\b\in S$. Now, we define $E(S)=E(S,R)$ as the subgroup generated by
all elementary root unipotent subgroups $X_{\a}$, $\a\in S$:
$$ E(S,R)=\langle x_{\a}(\xi),\quad \a\in S,\quad \xi\in R\rangle. $$
\noindent
In this notation, $U$ and $U^{-}$ coincide with $E(\Phi^{+},R)$ and
$E(\Phi^{-},R)$, respectively. The groups $E(S,R)$ are particularly
important in the case where $S$ is a {\it special\/} =
{\it unipotent\/} set of roots, in other words, where
$S\cap(-S)=\varnothing$. In this case $E(S,R)$ coincides with the
{\it product\/} of root subgroups $X_{\a}$, $\a\in S$, in some/any
fixed order.
\par
Let again $S\subseteq\Phi$ be a closed set of roots. Then $S$
can be decomposed into a disjoint union of its
{\it reductive\/} = {\it symmetric\/} part $S^{r}$,
consisting of those $\a\in S$, for which $-\a\in S$, and its
{\it unipotent\/} part $S^{u}$, consisting of those $\a\in S$, for
which $-\a\not\in S$. The set $S^{r}$ is a closed root subsystem,
whereas the set $S^{u}$ is special. Moreover, $S^{u}$
is an {\it ideal\/} of $S$, in other words, if $\a\in S$,
$\b\in S^{u}$ and $\a+\b\in\Phi$, then $\a+\b\in S^{u}$.
{\it Levi decomposition\/} asserts that the group $E(S,R)$ decomposes
into semidirect product $E(S,R)=E(S^r,R)\rightthreetimes E(S^u,R)$
of its {\it Levi subgroup\/} $E(S^{r},R)$ and its
{\it unipotent radical\/}~$E(S^{u},R)$.

\section*{\S6. Elementary parabolic subgroups}

\noindent
The main role in the proof of Theorem 1 is played by Levi
decomposition for elementary parabolic subgroups. Denote by
$m_k(\a)$ the coefficient of $\a_k$ in the expansion of $\a$
with respect to the fundamental roots:
$$ \a=\sum m_k(\a)\a_k,\quad 1\le k\le l. $$
\par
Now, fix an $r=1,\ldots,l$ -- in fact, in the reduction to smaller
rank it suffices to employ only terminal parabolic subgroups,
even only the ones corresponding to the first and the last
fundamental roots, $r=1,l$. Denote by
$$ S=S_r=\big\{\a\in\Phi,\ m_r(\a)\geq 0\big\} $$
\noindent
the $r$-th standard parabolic subset in $\Phi$. As usual,
the reductive part $\Delta=\Delta_r$ and the special part
$\Sigma=\Sigma_r$ of the set $S=S_r$ are defined as
$$ \Delta=\big\{\alpha\in\Phi,\ m_r(\alpha) = 0\big\},\quad
\Sigma=\big\{\alpha\in\Phi,\ m_r(\alpha) > 0\big\}. $$
\noindent
The opposite parabolic subset and its special part are defined
similarly
$$ S^-=S^-_r=\big\{\alpha\in\Phi,\ m_r(\alpha)\leq 0\big\},\quad
\Sigma^-=\big\{\alpha\in\Phi,\ m_r(\alpha)<0\big\}. $$
\noindent
Obviously, the reductive part $S^-_r$ equals $\Delta$.
\par
Denote by $P_r$ the {\it elementary\/} maximal parabolic
subgroup of the elementary group $E(\Phi,R)$. By definition,
$$ P_r=E(S_r,R)=\big\langle x_\alpha(\xi),\ \alpha\in S_r,
\ \xi\in R \big\rangle. $$
\noindent
Now Levi decomposition asserts that the group $P_r$ can be represented
as the semidirect product
$$ P_r=L_r\rightthreetimes U_r=E(\Delta,R)\rightthreetimes E(\Sigma,R) $$
\noindent
of the elementary Levi subgroup $L_r=E(\Delta,R)$ and the unipotent
radical $U_r=E(\Sigma,R)$. Recall that
$$ L_r=E(\Delta,R)=\big\langle x_\alpha(\xi),\quad \alpha\in\Delta,
\quad \xi\in R \big\rangle, $$
\noindent
Whereas
$$ U_r=E(\Sigma,R)=
\big\langle x_\alpha(\xi),\ \alpha\in\Sigma,\ \xi\in R\big\rangle. $$
\noindent
A similar decomposition holds for the opposite parabolic subgroup
$P_r^-$, whereby the Levi subgroup is the same as for $P_r$,
but the unipotent radical $U_r$ is replaced by the opposite unipotent
radical $U_r^-=E(-\Sigma,R)$
\par
As a matter of fact, we use Levi decomposition in the following
form. It will be convenient to slightly change the notation
and write $U(\Sigma,R)=E(\Sigma,R)$ and $U^-(\Sigma,R)=E(-\Sigma,R)$.

\begin{proclaim}{Lemma 3}
The group\/ $\big\langle U^{\sigma}(\Delta,R),U^\rho(\Sigma,R)\big\rangle$,
where\/ $\sigma,\rho=\pm 1$, is the semidirect product of its
normal subgroup\/ $U^\rho(\Sigma,R)$ and the complementary subgroup\/
$U^{\sigma}(\Delta,R)$.
\end{proclaim}

In other words, it is asserted here that the subgroup $U^{\pm}(\Delta,R)$
normalised each of the groups $U^{\pm}(\Sigma,R)$, so that, in
particular, one has the following four equalities for products
$$ U^{\pm}(\Delta,R)U^{\pm}(\Sigma,R)=U^{\pm}(\Sigma,R)U^{\pm}(\Delta,R), $$
\noindent
and, furthermore, the following four obvious equalities for
intersections hold:
$$ U^{\pm}(\Delta,R)\cap U^{\pm}(\Sigma,R)=1. $$
\par
In particular, one has the following decompositions:
$$ U(\Phi,R)=U(\Delta,R)\rightthreetimes U(\Sigma,R),
\quad
U^-(\Phi,R)=U^-(\Delta,R)\rightthreetimes U^-(\Sigma,R). $$

\section*{\S7. Reduction to groups of smaller rank}

The following result is a minor elaboration of Proposition 1
from the paper by Oleg Tavgen \cite{tavgen90b}. Tavgen states a slightly
weaker, but more general result, in terms of existence of
factorisations for {\it all\/} irreducible root systems of
a certain rank\footnote{Observe, that one should read the
first equality in the statement of Proposition 1 in \cite{tavgen90b} as
$\rk({}^{\sigma}\Phi_0)=m$.}. On the other hand, he considers
also twisted groups, apart from those of type ${}^2\!\!\A_{2l}$.

\begin{proclaim}{Theorem 3}
Let\/ $\Phi$ be a reduced irreducible root system of rank
$l\ge 2$, and\/ $R$ be a commutative ring. Suppose that
for subsystems\/ $\Delta=\Delta_1,\Delta_l$ the elementary
Chevalley group\/ $E(\Delta,R)$ admits unitriangular
factorisation
$$ E(\Delta,R)=(U(\Delta,R)U^-(\Delta,R))^L. $$
\noindent
Then the elementary Chevalley group\/ $E(\Phi,R)$ admits
unitriangular factorisation
$$ E(\Phi,R)=(U(\Phi,R)U^-(\Phi,R))^L. $$
\noindent
of the same length\/ $2L$.
\end{proclaim}

Clearly, Theorem~1 immediately follows from Lemma 1 and
Theorem 3, so that it only remains to prove Theorem~3.
\par
The leading idea of Tavgen's proof is so general and beautiful
that it works in many other similar contexts. It relies
on the fact that for systems of rank $\ge 2$ every fundamental
root falls into the subsystem of smaller rank obtained by
dropping either the first or the last fundamental root.
Eiichi Abe and Kazuo Suzuki \cite{abe69} and \cite{AS} used the same
argument in their description of normal subgroups,
to extract root unipotents. A similar consideration, in
conjunction with a general position argument, was used by
Vladimir Chernousov, Erich Ellers, and Nikolai Gordeev
in their simplified proof of the Gauss decomposition with
prescribed semisimple part \cite{CEG}.
\par
Let us reproduce the details of the argument. By definition
$$ Y=(U(\Phi,R)U^-(\Phi,R))^L $$
\noindent
is a {\it subset\/} in $E(\Phi,R)$. Usually, the easiest way
to prove that a subset $Y\subseteq G$ coincides with the whole
group $G$ consists in the following.

\begin{proclaim}{Lemma 4}
Assume that\/ $Y\subseteq G$, $Y\neq \varnothing$,
and\/ $X\subseteq G$ be a symmetric generating set. If\/
$XY\subseteq Y$, then\/ $Y=G$.
 \end{proclaim}

\begin{proof}[Proof of Theorem $3$]
By Lemma 2 the group $G$ is generated by the fundamental
root elements
$$ X=\big\{x_{\a}(\xi)\mid \a\in\pm\Pi,\ \xi\in R\big\}. $$
\noindent
Thus, by Lemma 4 is suffices to prove that $XY\subseteq Y$.
\par
Let us fix a fundamental root unipotent $x_{\a}(\xi)$.
Since $\rk(\Phi)\ge 2$, the root $\a$ belongs to at least one
of the subsystems $\Delta=\Delta_r$, where $r=1$ or $r=l$,
generated by all fundamental roots, except for the first or the
last one, respectively. Set $\Sigma=\Sigma_r$ and express
$U^{\pm}(\Phi,R)$ вin the form
$$ U(\Phi,R)=U(\Delta,R)U(\Sigma,R),\quad
U^-(\Phi,R)=U^-(\Delta,R)U^-(\Sigma,R). $$
\par
Using Lemma 3 we see that
$$ Y=(U(\Delta,R)U^-(\Delta,R))^L(U(\Sigma,R)U^-(\Sigma,R))^L. $$
\noindent
Since $\a\in\Delta$, one has $x_{\a}(\xi)\in E(\Delta,R)$, so that
the inclusion $x_{\a}(\xi)Y\subseteq Y$ immediately follows from
the assumption.
\end{proof}

\section*{\S8. The use of Bruhat and Gauss decomposition}

Observe that, as opposed to the proofs in \cite{BNP}, the proof
in the paper by Liebeck and Pyber \cite{LP} is natural and
transparent, and is based on the Bruhat decompsoition
$g=udn_wv$, where $u,v\in U$, $d\in T$, $w\in W$.
Excessive number of factors in their result is due to the fact
that they decompose $d$ and $n_w$ separately, and, on top of that,
the way they do it is far from being optimal. For instance,
to decompose an element $h_{\a}(\e)$ they rely directly on the
definition $h_{\a}(\e)=w_{\a}(\e)w_{\a}(-1)$, rather than
perform the actual calculation in $\SL(2,R)$. Substituting the
expression of $w_{\a}(\e)$ in terms of root unipotents, in the
definition of $h_{\a}(\e)$ we express $h_{\a}(\e)$ as a product
of 5 root unipotents as follows:
$$ \begin{pmatrix} \e&0\\ 0&\e^{-1}\\  \end{pmatrix}=
 \begin{pmatrix} 1&\e\\ 0&1\\  \end{pmatrix}
 \begin{pmatrix} 1&0\\ -\e^{-1}&1\\  \end{pmatrix}
 \begin{pmatrix} 1&\e-1\\ 0&1\\  \end{pmatrix}
 \begin{pmatrix} 1&0\\ 1&1\\  \end{pmatrix}
 \begin{pmatrix} 1&-1\\ 0&1\\  \end{pmatrix}.  $$
\noindent Factually, it is well known that $h_{\a}(\e)$ is already a
product of 4 root unipotents:
$$ \begin{pmatrix} \e&0\\ 0&\e^{-1}\\  \end{pmatrix}=
 \begin{pmatrix} 1&-1\\ 0&1\\  \end{pmatrix}
 \begin{pmatrix} 1&0\\ 1-\e&1\\  \end{pmatrix}
 \begin{pmatrix} 1&\e^{-1}\\ 0&1\\  \end{pmatrix}
 \begin{pmatrix} 1&0\\ \e(\e-1)&1\\  \end{pmatrix}.  $$
\noindent
As a matter of fact, the same idea allows to {it immediately\/}
obtain for Chevalley groups over semilocal rings unitriangular
fcatorisation of the same length 5, as in \cite{BNP}. The proof is
based on Gauss decomposition and the following toy version of
Theorem~1.
\par
Recall, that the following analogue of Gauss decomposition was
established by Eiichi Abe and Kazuo Suzuki \cite{abe69,AS} and by
Michael Stein \cite{stein}\footnote{In terms of computational linear
algebra this is the ULU-decomposition, whereas Bruhat decomposition
is the UPU-decomposition, where P should be interpreted as the
first letter of the word `permutation'. Though, Stein himself
called this decomposition a {\it Bruhat type decomposition\/}.
Erich Ellers and Nikolai Gordeev usually intend by Gauss
decomposition the LPU-decomposition, which is normally called
{\it Birkhoff decomposition\/}. At the same time, in the works
on computational linear algebra the name `Gauss decomposition'
often refers either to the PLU-decomposition, or the
LUP-decomposition.}.

\begin{proclaim}{Lemma~5}
Let\/ $R$ be a semilocal ring. Then for a simply connected
Chevalley group one has the following decomposition
$$ G(\Phi,R)=U(\Phi,R)T(\Phi,R)U^-(\Phi,R)U(\Phi,R). $$
\end{proclaim}

As we know, Gauss decomposition for the group $\SL(n,R)$ holds under
a weaker assumption $\sr(R)=1$, see, {\it for instance\/}, \cite{VW2}.
First author \cite{vavilov} noticed that in fact condition $\sr(R)=1$ 
is {\it
necessary\/} for a Chevalley group $G(\Phi,R)$ to admit Gauss
decomposition, for {\it all\/} types. For the linear case this
obvious circumstance was rediscovered some 20+ later in \cite{NDS,CC}. On
the other hand, for all types except $\Phi=\A_l,\C_l$ all known
sufficient conditions are somewhat stronger, say, something like
$\asr(R)=1$. Thus, there is still some gap between the necessary and
sufficient conditions.
\par
The following result is a generalisation of Lemma 18 of \cite{DV88},
where a slightly weaker fact is proven for the case $\Phi=\A_l$.

\begin{proclaim}{Theorem 4}
Let\/ $\Phi$ be a reduced irreducible root system and\/ $R$
be an arbitrary commutative ring. Then one has the following
inclusion
$$ N(\Phi,R)\subseteq U(\Phi,R)U^-(\Phi,R)U(\Phi,R)U^-(\Phi,R). $$
\end{proclaim}

\begin{proclaim}{Corollary}
Let\/ $\Phi$ be a reduced irreducible root system and\/ $R$
be a commutative semilocal ring. Then the simply connected
Chevalley group\/ $G(\Phi,R)$ admits unitriangular factorisation
$$ G(\Phi,R)=U(\Phi,R)U^-(\Phi,R)U(\Phi,R)U^-(\Phi,R)U(\Phi,R) $$
\noindent
of length\/ $5$.
\end{proclaim}

\begin{proof}[Proof]
Indeed,
$$ G=UTU^-U\le U(UU^-UU^-)U^-U=UU^-UU^-U. $$
\end{proof}

Since this result is both less general in terms of condition
on the ground ring and weaker in terms of the resulting length
than Theorem~1, we do not present the proof of Theorem~4
in general. However, for the purpose of illustration, let us
reproduce its proof in the linear case. Recall that in this
case the group $N(\A_{n-1},R)$ coincides with the group
$N=\SN(n,R)$ of monomial matrices with determinant 1.

\begin{proof}[Proof of Theorem $4$ for $\Phi=\A_{n-1}$]
Let $g=(g_{ij})\in\SN(n,R)$. Let us argue by induction on $n$.
In the case $n=1$ there is nothing to prove. Thus, let $n\ge 2$.

\medskip\noindent
{\bf Case 1.} First, let $g_{nn}=0$. Then there exists a unique
$1\le r\le n-1$ such that $a=g_{rn}\neq 0$ and a unique $1\le s\le
n-1$ such that $b=g_{ns}\neq 0$, all other entries in the $s$-th and
the $n$-th columns are equal to 0. Since $g$ is invertible,
automatically $a,b\in R^*$. The matrix $gt_{sn}(b^{-1})$ differs
from $g$ only in the position $(n,n)$, where now we have 1 instead
of 0. Consecutively multiplying the resulting matrix on the right by
$t_{ns}(-b)$ and then by $t_{sn}(b^{-1})$, we get the matrix $h$,
which differs from $g$ only at the intersection of the $r$-th and
the $n$-the rows with the $s$-th and the $n$-th columns, where now
instead of $ \begin{pmatrix} 0&a\\ b&0\\  \end{pmatrix}$ one has $
\begin{pmatrix} -ab&0\\ 0&1\\  \end{pmatrix}$. Observe, that the
determinant of the leading submatrix of order $n-1$ of the matrix
$h$ equals 1, and thus we can apply induction hypothesis and obtain
for that last matrix the desired factorisation in the group
$\SL(n-1,R)$. This factorisation does not affect the last row and
the last column. Thus, all factors of the above factorisation of the
submatrix $h$ lie in $L_{n-1}$. It only remains to notice that
$t_{sn}(b^{-1})\in U_{n-1}$, and $t_{ns}(-b)\in U_{n-1}^-$, and then
to invoke Lemma~3.

\medskip\noindent
{\bf Case 2.} Now, let $b=g_{nn}\neq 0$. Take arbitrary
$1\le r,s\le n-1$ for which $a=g_{rs}\neq 0$. Again,
automatically $a,b\in R^*$. As in the previous case, let us
concentrate on the $r$-th and the $n$-th rows and the
$s$-th and the $n$-th columns. Since there are no further
non-zero entries in these rows and columns, any additions
between them do not change other entries of the matrix, and
only affect the submatrix at the intersection of the $r$-th
and the $n$-th rows with the $s$-th and the $n$-th columns.
Now, multiplying $g$ by
$t_{ns}(b^{-1})t_{sn}(1-b)t_{ns}(-1)t_{ns}(-b^{-1}(1-b))$
on the right, we obtain the matrix $h$, where this submatrix,
which was initially equal to
$\begin{pmatrix} a&0\\ 0&b\\  \end{pmatrix}$, will be replaced
by $\begin{pmatrix} ab&0\\ 0&1\\  \end{pmatrix}$. At this point
the proof can be finished in exactly the same way as in the
previous case.
\end{proof}

\section*{\S9. Dedekind rings of arithmetic type}

\def\p{\mathfrak p}
Let $K$ be a global field, i.~e.\ either a finite algebraic
extension of the field $\Rat$, or a field of algebraic functions
in one variable over a finite field of constants $\GF{q}$.
Further, let $S$ be a finite set of non-equivalent valuations of
$K$, nonempty in the functional case and containing all Archimedian
valuations in the number case. For a non-Archimedian valuation
$\p$ of the field $K$ we denote by $v_{\p}$ the corresponding exponent.
\par
As usual, $R={\mathcal O}_S$ denotes the ring consisting of all
$x\in K$ such that $v_{\p}(x)\ge 0$ for all valuations $\p$ of the
field $K$, which do not belong to $S$. The ring ${\mathcal O}_S$
is called the Dedekind ring of arithmetic type determined by
the set of valuations $S$ of the field $K$, or, otherwise, a
Hasse domain, see, for instance, \cite{BMS}. We will be mostly 
interested in the case $|S|\ge 2$, where by Dirichlet unit theorem 
the ring ${\mathcal O}_S$ has a unit of infinite order.
\par
As another immediate corollary of reduction to smaller ranks, of
the positive solution for the congruence subgroup problem 
\cite{BMS,matsumoto,serre}
and of the Cooke---Weinberger paper \cite{CW}, one can state the
following result. Observe, that Cooke---Weinberger calculations
depend on the infinity of primes in arithmetic progressions,
subject to additional multiplicative restrictions, and thus
depend upon GRH = Generalised Riemann's Hypothesis.

\begin{proclaim}{Theorem 5}
Let\/ $\Phi$ be a reduced irreducible root system and\/
$R={\mathcal O}_S$ be a Dedekind ring of arithmetic type
with an infinite multiplicative group. Under assumption of the
Generalised Riemann's Hypothesis the simply connected Chevalley
group\/ $G(\Phi,R)$ admits unitriangular factorisation
$$ G(\Phi,R)=(U(\Phi,R)U^-(\Phi,R))^4U(\Phi,R) $$
\noindent
of length\/ $9$.
\end{proclaim}

One can prove similar results also in the absence of Generalised
Riemann's Hypothesis. However, in that case the known length
estimates would be much worse, and depend on the number of classes
of the ring $R$, the number of prime divisors of its discriminant,
or something of the sort. Nevertheless, we are convinced that 9 here
could be actually replaced by 5 or 6.

\begin{proof}[Proof of Theorem~$2$]
Theorem~3 implies that one only has to prove the existence
of a unitriangular factorisation of length 6 for the group
$\SL\Big(2,\Int\Big[\displaystyle{1\over p}\Big]\Big)$.
In fact, we will prove the following slightly more precise result,
where
$U=U\Big(2,\Int\Big[\displaystyle{1\over p}\Big]\Big)$ and
$U^-=U^-\Big(2,\Int\Big[\displaystyle{1\over p}\Big]\Big)$,
respectively.
\end{proof}

\begin{proclaim}{Lemma 6}
Let\/ $p\in\Int$ be a rational prime. Then under assumption
of the Generalised Riemann's Hypothesis one has
$$ \SL\Big(2,\Int\Big[\displaystyle{1\over p}\Big]\Big)=
U^-UU^-UU^-\enskip\textstyle{\bigcup}\enskip UU^-UU^-U. $$
\end{proclaim}

Our proof critically depends on the recent results on Artin's
conjecture, specifically, on the following result which was
first stated in the papers by Pieter Moree \cite{moree99,moree08}.
A complete proof can be found in \cite{LMS}, Corollary~5.4.  

\begin{proclaim}{Lemma 7}
Let\/ $a\in\Int$, $a\neq 0,1,-1$, is a square-free integer
and $c,d\in\Int$, $c\perp d$, be coprime integers. Then
under assumption of the Generalised Riemann's Hypothesis
density of the set of primes $q$ in the residue class
$c\pmod{d}$, for which $q$ is a primitive root $\mod d$,
exists and is positive, with the only exception of the case,
where the discriminant of $\Rat(\sqrt{a})$ divides $d$, and
every prime in the residue class $c\pmod{d}$ completely
decomposes in $\Rat(\sqrt{a})$.
\end{proclaim}

Let us state the following immediate corollary of this result.

\begin{proclaim}{Corollary}
Let\/ $p\in\Int$ be a rational prime and $c\perp d$ be two coprime
integers such that $p\perp d$. Then under assumption of the
Generalised Riemann's Hypothesis there are infinitely many primes
$q$ in the residue class $c\pmod{d}$, for which $p$ is a primitive
root modulo $q$.
\end{proclaim}
Now, we are all set to prove Lemma~6, and thus also Theorem~2.

\begin{proof}[Proof]
Let $g=\begin{pmatrix} x&y\\ z&w\\  \end{pmatrix}\in
\SL\Big(2,\Int\Big[\displaystyle{1\over p}\Big]\Big)$.
Since the case, where at least one of the matrix entries
$x,y,z,w$ equals 0 was already considered in the above
proof of Theorem~4 for the linear case, in the sequel we
can assume that $xyzw\neq 0$. Now, set
$$ g= \begin{pmatrix} x&y\\ z&w\\  \end{pmatrix}=
 \begin{pmatrix} p^{\alpha}a & p^{\beta}b\\ \ast&\ast\\  \end{pmatrix}\in
\SL\Big(2,\Int\Big[\displaystyle{1\over p}\Big]\Big), $$
\noindent
where $a,b\in\Int$ are not divisible by $p$, and $\alpha,\beta\in\Int$.

\noindent
{\bf Case 1: $\alpha\geq\beta$}. By corollary of Lemma 7, there
are infinitely many primes $q$ of the form $p^{\alpha-\beta}a+bk$
such that $p$ is a primitive root modulo $q$. Then
$$ gt_{21}(k)= \begin{pmatrix}
p^{\beta}q & p^{\beta}b\\ \ast&\ast\\  \end{pmatrix} .  $$
\noindent
Since $p$ is a primitive root $\mod{q}$, there exists a
$u\geq 1$ such that $p^u\equiv b\pmod{q}$. Let, say,
$p^u=b+lq$. In this case
$$ gt_{21}(k)t_{12}(l)= \begin{pmatrix}
p^{\beta}q&p^{\beta+u}\\ \ast&\ast\\  \end{pmatrix}. $$
\noindent
Then for $\theta=(1-p^{\alpha}q)/p^{\beta+u}$, we get
$$ gt_{21}(k)t_{12}(l)t_{21}(\theta)=
\begin{pmatrix} 1&p^{\beta+u}\\ \ast&\ast\\  \end{pmatrix} $$
\noindent
and, finally,
$$ gt_{21}(k)t_{12}(l)t_{21}(\theta)t_{12}(-p^{\beta+u})=
\begin{pmatrix} 1&0\\ \ast&1\\  \end{pmatrix}. $$
\noindent
Thus, in this case $g\in U^-UU^-UU^-$.

\medskip\noindent
{\bf Case 2: $\alpha<\beta$.} In this case exactly the same
argument shows that $g\in UU^-UU^-U$.
\end{proof}

\section*{\S10. Some related unsolved problems}

Let us list some further problems which are intimately
related to unitriangular factorisations. Our Theorem~1
allows to slightly improve known bounds in some of them. 

\smallskip
$\bullet$ Let us mention another problem similar to the one
discussed in the present paper. Estimate the width of the
elementary group $E(\Phi,R)$ with respect to {\it all\/}
unipotent elements. From Gauss decomposition with prescribed
semisimple part, as obtained by Erich Ellers and Nikolai
Gordeev, it follows that every noncental element of a
Chevalley group $G(\Phi,K)$ over a field $K$ is a product
of two unipotent elements \cite{EG}.

\smallskip
Theorem 1 immediately implies the following result.

\begin{proclaim}{Corollary 1}
Let\/ $\Phi$ be a reduced irreducible root system and\/ $R$
be a commutative ring such that $\sr(R)=1$. Then every
element of the simply connected Chevalley group $G(\Phi,R)$
is a product of three unipotent elements.
\end{proclaim}

For rings of dimension $\ge 1$, the situation is much more delicate,
see, for example, the paper by Fritz Grunewald, Jens Mennicke and
Leonid Vaserstein \cite{GMV}.

\smallskip
$\bullet$ Another related problem, which was recently    
considered in the paper by Martin Liebeck, Nikolay Nikolov and Aner
Shalev \cite{LNS}, is the width of a Chevalley group in fundamental
$\SL(2,R)$'s. For Chevalley groups over a [finite] field that paper
gives the estimate $5|\Phi^+|$, which follows from \cite{BNP}. Our
Theorem~1 immediately implies a better bound in a much more general
situation.

\begin{proclaim}{Corollary 2}
Let\/ $\Phi$ be a reduced irreducible root system and\/ $R$
be a commutative ring such that $\sr(R)=1$. Then the simply
connected Chevalley group $G(\Phi,R)$ can be written as the
product of $4|\Phi^+|$ copies of the fundamental $\SL(2,R)$.
\end{proclaim}

As a matter of fact, even this estimate is far from the real one.
Bruhat decomposition immediately implies the bound $3|\Phi^+|$, and
somewhat more precise calculations allow to get the estimate
$2|\Phi^+|$, for all Bezout rings. This is done in the paper by the
first author and Evdokim Kovach, currently under way.

\smallskip
$\bullet$ A large number of papers are dedicated to the width
of Chevalley groups in commutators, see in particular
\cite{FV88,VW2,you,YZ,AVG,EG,LOST10,SS,SV}. For Chevalley groups 
over a ring of stable rank 1 our Theorem~1 implies the following
obvious estimate, 
\begin{proclaim}{Corollary 3}
Let\/ $\Phi$ be a reduced irreducible root system and\/ $R$
be a commutative ring such that $\sr(R)=1$. Then the commutator 
width of the elementary Chevalley groups $E(\Phi,R)$ does not
exceed $5$.
\end{proclaim}
\begin{proof}[Proof]
For any $g\in E(\Phi,R)$ there exist $u,v\in U(\Phi,R)$ and
$x,y\in U^-(\Phi,R)$ such that
$$ g=uxvy=\big[{}^ux,{}^uv]\cdot uv\cdot xy. $$
\noindent
Since $uv\in U(\Phi,R)$ and $xy\in U^-(\Phi,R)$, the claim follows.
\end{proof}
We believe that similarly to \cite{you,AVY} the commutator width 
of elementary Chevalley groups over a ring of stable rank 1 does 
not exceed 2 or 3.

\smallskip
In conclusion, let us state some further unsolved problems in
the field.

\begin{proclaim}{Problem 1}
Find for a Chevalley group of rank\/ $\ge 2$ the minimal\/ $L$
such that
$$ G(\Phi,\Int)=(U(\Phi,\Int)U^-(\Phi,\Int))^L. $$
\end{proclaim}

In fact, the precise estimate is not even known in the following
special case. We believe that the estimate by Thomas Laffey
\cite{laffey-a,laffey-b}, which follows from \cite{DV88}, is 
grossly exaggerated.

\begin{proclaim}{Problem 2}
Find the minimal\/ $n$, starting from which one has the
unitriangluar factorisation
$$ \SL(n,\Int)=(U(n,\Int)U^-(n,\Int))^3. $$
\end{proclaim}

The following question is of obvious relevance, in connection
with Lemma 6.

\begin{proclaim}{Problem 3}
Is it true that $U^-UU^-UU^-=UU^-UU^-U${\rm?}
\end{proclaim}

Most probably, in general this is not the case, but it would
be very interesting to see an explicit counter-example.
\par
Generally speaking, no unitriangular factorisations exist over
rings of dimension $>1$, and no such factorisations can possibly
exist. Their role is taken by parabolic factorisations, see
\cite{VS10,VS11,SiV}, and the references therein. These factorisations
allow to efficiently reduce problems regarding the group itself,
to similar problems for groups of smaller ranks.

\begin{proclaim}{Problem 4}
Estimate the width of Chevalley groups over commutative rings of
small stable ranks in terms of classical subgroups/subgroups of
type\/ $\A_l$.
\end{proclaim}

Even for fields such estimates should be much better than the
ones obtained in \cite{nikolov,NP}.
\par
The following problem can be easily solved by the methods of the
present paper.

\begin{proclaim}{Problem 5}
Calculate the width of the elementary Chevalley group\/
$E(\Phi,R)$ over a semilocal ring\/ $R$ in terms of unipotent
radicals\/ $U_P$ and $U_P^-$ of two opposite parabolic subgroups.
\end{proclaim}

Let us mention yet another direction in which it would be
natural to generalise the results of the present paper.
In \cite{PS} Victor Petrov and Anastasia Stavrova constructed the
elementary subgroup $E(R)$ in an isotropic reductive group
$G(R)$ and, for groups of relative rank $\ge 2$, proved
normality of $E(R)$ in $G(R)$. See also the recent paper \cite{LuS}
by Alexander Luzgarev and Anastasia Stavrova, where it is
proven that $E(R)$ is perfect, with the known exceptions for
Chevalley groups of rank 2.
\par
Despite the apparent similarity of the statements, the following
problem is not at all trivial, and, to the best of our knowledge,
presently there are no non-trivial estimates at all.

\begin{proclaim}{Problem 6}
Calculate the width of the elementary subgroup\/ $E(R)$ of an
isotropic reductive group\/ $G(R)$ over a semilocal ring\/ $R$,
in terms of unipotent radicals\/ $U_P$ and\/ $U_P^-$ of two
opposite parabolic subgroups.
\end{proclaim}

Even less than that is known in the arithmetic context.

\begin{proclaim}{Problem 7}
Prove that the elementary subgroup\/ $E(R)$ of an isotropic
reductive group\/ $G(R)$ of relative rank\/ $\ge 2$ has bounded
width with respect to the unipotent radicals\/ $U_P$ and\/ $U_P^-$
of two opposite parabolic subgroups, in the case where $R={\mathcal
O}_S$ is a Dedekind ring of arithmetic type.
\end{proclaim}

The only result we are aware of, regarding this problem, is the
paper by Igor Erovenko and Andrei Rapinchuk \cite{ER}, which addresses
the case of orthogonal groups, corresponding to a form of Witt
index $\ge 2$.
\par
\medskip
The authors thank Pieter Moree for references related to Artin's
conjecture, and Dave Witte Morris for extremely useful discussions
of different approaches towards the proof of bounded generation.

\smallskip

Vavilov N.~A., Smolensky A.~V., Sury B.
Unitriangular factorisations of Chevalley groups.

Lately, the following problem attracted a lot of attention
in various contexts: find the shortest factorisation
$G=UU^-UU^-\ldots U^{\pm}$ of a Chevalley group
$G=G(\Phi,R)$ in terms of the unipotent radical $U=U(\Phi,R)$
of the standard Borel subgroup $B=B(\Phi,R)$ and the unipotent
radical $U^-=U^-(\Phi,R)$ of the opposite
Borel subgroup $B^-=B^-(\Phi,R)$. So far, the record over a
finite field was established in a 2010 paper by Babai,
Nikolov, and Pyber, where they prove that a group
of Lie type admits unitriangular factorisation $G=UU^-UU^-U$
of length 5. Their proof invokes deep analytic and
combinatorial tools. In the present paper we notice that
from the work of Bass and Tavgen one immediately gets a much
more general result, asserting that over any ring of stable
rank 1 one has unitriangular factorisation $G=UU^-UU^-$
of length 4. Moreover, we give a detailed survey of
triangular factorisations, prove some related results,
discuss prospects of generalisation to other classes of rings,
and state several unsolved problems. Another main result
of the present paper asserts that, in the assumption of
the Generalised Riemann's Hypothesis, Chevalley groups
over the ring $\Int\Big[\displaystyle{1\over p}\Big]$ admit
unitriangular factorisation $G=UU^-UU^-UU^-$ of length 6.
Otherwise, the best length estimate for Hasse domains with
infinite multiplicative groups that follows from the work
of Cooke and Weinberger, gives 9 factors.


\begin{thebibliography}{99}

\bibitem{BMS} H.~Bass, J.~Milnor, J.-P.~Serre,
{\it Solution of the congruence subgroup problem for\/ $\SL_n$\/
{\rm ($n\ge 3$)} and\/ $\Sp_{2n}$\/ {\rm($n\ge 2$)}\/}.
Publ. Math. Inst. Hautes Etudes Sci. 33 (1967), 59--137.

\bibitem{vavilov84} N.~A.~Vavilov,
{\it Parabolic subgroups of Chevalley groups over a commutative
ring}. --- J. Soviet Math {\bf 116} (1984), 1848--1860.

\bibitem{VS10} N.~A.~Vavilov, S.~S.~Sinchuk,
{\it Dennis---Vaserstein type decomposition}.
--- Zapiski Nauchn. Semin. POMI {\bf375} (2010), 48--60
(in Russian, English translation pending).

\bibitem{VS11} N.~A.~Vavilov, S.~S.~Sinchuk,
{\it Parabolic factorisations of the split classical groups}.
--- Algebra and Analysis {\bf 23} (2011), No.~4, 1--30
(in Russian, English translation pending).

\bibitem{LuS} A.~Yu.~Luzgarev, A.~K.~Stavrova,
{\it Elementary subgroup of an isotropic reductive group is perfect}.
--- Algebra and Analysis {\bf23} (2011) (to appear).
(in Russian, English translation pending).

\bibitem{PS} V.~A.~Petrov, A.~K.~Stavrova,
{\it Elementary subgroups of isotropic reductive groups}.
--- St. Petersburg Math. J. {\bf 20}  (2008), no. 3, 160--188.

\bibitem{serre} J.-P.~Serre,
{\it Le probl\`eme des groupes de congruence pour\/ $\SL_2$}.
--- Ann. Math. {\bf 92} (1970), 489--527.

\bibitem{steinberg} R.~Steinberg, {\it Lectures on Chevalley groups},
Yale University, 1967.

\bibitem{tavgen90a}
O.~I.~Tavgen, {\it Finite width of arithmetic subgroups of Chevalley
groups of rank\/ $\ge 2$}. --- Soviet Math.\ Doklady
{\bf 41} (1990), no.1, 136--140.

\bibitem{tavgen90b}
O.~I.~Tavgen, {\it Bounded generation of Chevalley groups over rings
of\/ $S$-integer algebraic num\-bers\/}. --- Izv. Acad. Sci. USSR
{\bf 54} (1990), no.1, 97--122.

\bibitem{abe69} E.~Abe,
{\it Chevalley groups over local rings}.
 ---  T\^ohoku Math.  J. {\bf21} (1969), No.~3, 474--494.

\bibitem{AS} E.~Abe, K.~Suzuki,
{\it On normal subgroups of Chevalley groups over commutative rings}.
---  T\^ohoku Math. J. {\bf28} (1976), No.~1, 185--198.

\bibitem{AVY} F.~A.~Arlinghaus, L.~N.~Vaserstein, You Hong,
{\it Commutators in pseudo-orthogonal groups\/}. ---  J. Austral.
Math. Soc., Ser.~A {\bf 59} (1995), 353--365.

\bibitem{BNP} L.~Babai, N.~Nikolov, L.~Pyber,
 {\it Product growth and mixing in finite groups}.
--- In: 19th Annual ACM--SIAM Symposium on Discrete Algorithms,
ACM--SIAM (2008), pp.~248--257.

\bibitem{bass64} H.~Bass,
{\it $\K$-theory and stable algebra}.
---   Publ. Math. Inst. Hautes \'Etudes Sci. No.~22 (1964), 5--60.

\bibitem{CK83} D.~Carter, G.~Keller,
{\it Bounded elementary generation of\/ $\SL_n({ \mathcal O})$}.
---  Amer. J. Math. {\bf 105} (1983), 673--687.

\bibitem{CK84} D.~Carter, G.~Keller,
{\it Elementary expressions for unimodular matrices}. ---
Commun. Algebra {\bf12} (1984), 379--389.

\bibitem{CKP}  D.~Carter, G.~E.~Keller, E.~Paige,
 {\it Bounded expressions in\/ $\SL(2,{ \mathcal O})$}.
--- Preprint Univ. Virginia  (1983).

\bibitem{carter} R.~W.~Carter,
{\it  Simple groups of Lie type}.
Wiley, London et al., 1972.

\bibitem{CK} Chen Baoquan, A.~Kaufman,
 {\it {\rm 3D} volume rotation using shear transformations}.
---  Graph. Mo\-dels {\bf62} (2000), 308--322.

\bibitem{CC} Chen Huanyin, Chen Miaosen,
 {\it On products of three triangular matrices over
associative rings}.
---  Linear Algebra Applic. {\bf387} (2004), 297--311.

\bibitem{CEG} V.~Chernousov, E.~Ellers, N.~Gordeev,
 {\it Gauss decomposition with prescribed semisimple part:
short proof}.
---  J. Algebra {\bf229} (2000), 314--332.

\bibitem{cohn} P.~M.~Cohn,
{\it On the structure of the\/ $\GL_2$ of a ring}.
---  Publ. Math. Inst. Hautes \'Etudes Sci.  No.~30 (1967), 5--53.

\bibitem{CW} G.~Cooke, P.~J.~Weinberger,
{\it On the construction of division chains in algebraic number rings,
with applications to\/ $\SL_2$}.
 ---  Commun. Algebra  {\bf3} (1975), 481--524.

\bibitem{DV88} R.~K.~Dennis, L.~N.~Vaserstein,
{\it On a question of M.~Newman on the number of commutators}.
--- J.~Algebra {\bf118} (1988), 150--161.

\bibitem{EG}  E.~Ellers, N.~Gordeev,
 {\it On the conjectures of J.~Thompson and O.~Ore}.
---  Trans. Amer. Math. Soc. {\bf350} (1998), 3657--3671.

\bibitem{ER} C.~K.~Fong, A.~R.~Sourour,
{\it The group generated by unipotent operators}.
--- Proc. Amer. Math. Soc. {\bf 97} (1986), No.~3, 453--458.

\bibitem{FS} I.~V.~Erovenko, A.~S.~Rapinchuk,
{\it Bounded generation of some\/ $S$-arithmetic orthogonal groups}.
---  C.~R.~Acad.~Sci. {\bf333} (2001), No.~5, 395--398.

\bibitem{GMV}  F.~J.~Grunewald, J.~Mennicke, L.~N.~Vaserstein,
 {\it On the groups\/ $\SL_2(\Int[x])$ and\/ $\SL_2(K[x,y])$}.
---  Israel J. Math. {\bf86} (1994), Nos.~1--3, 157--193.

\bibitem{GM}  R.~M.~Guralnick, G.~Malle,
 {\it Products of conjugacy classes and fixed point spaces}.
  {\tt arXiv: 1005.3756}.

\bibitem{hao} Hao Pengwei,
 {\it Customizable triangular factorizations of matrices}.
---  Linear Algebra Applic. {\bf382} (2004), 135--154.

\bibitem{vdK82} W.~van der Kallen,
{\it $\SL_3(\Co[x])$ does not have bounded word length}.
--- Springer Lect. Notes Math. {\bf966} (1982), 357--361.

\bibitem{laffey-a} T.~J.~Laffey,
 {\it Expressing unipotent matrices over rings as products of involutions}.
--- Preprint Univ. Dublin (2010).

\bibitem{laffey-b}  T.~J.~Laffey,
 {\it Lectures on integer matrices}.
---  Beijing (2010), 1--38.

\bibitem{LHW}  Lei Yang, Hao Pengwei, Wu Dapeng,
 {\it Stabilization and optimization of PLUS factorization
and its application to image coding}.
---  J. Visual Communication \& Image Representation {\bf22}
(2011), No.~1.

\bibitem{LS} M.~Larsen, A.~Shalev,
 {\it Word maps and Waring type problems}.
---  J. Amer. Math. Soc. {\bf22} (2009), 437--466.

\bibitem{LMS} H.~W.~Lenstra (jr.), P.~Moree, P.~Stevenhagen,
{\it Character sums for primitive root densities}.
---  (2011) (to appear).

\bibitem{LS96}  M.~Liebeck, A.~Shalev,
 {\it Classical groups, probabilistic methods, and
the $(2,3)$-generation problem}.
---  Ann. Math. {\bf144} (1996), No.~1, 77-125.

\bibitem{LS01}  M.~Liebeck, A.~Shalev,
 {\it Diameteres of finite simple groups: sharp bounds and
applications}.
---  Ann. Math. {\bf154} (2001), 383--406.

\bibitem{LNS}  M.~Liebeck, N.~Nikolov, A.~Shalev,
 {\it Groups of Lie type as products of $\SL_2$ subgroups}.
 --- J. Algebra {\bf 326} (2011), 201--207.

\bibitem{LOST10}  M.~Liebeck, E.~A.~O'Brien, A.~Shalev, Pham Huu Tiep,
 {\it The Ore conjecture}.
---  J. Europ. Math. Soc. {\bf12} (2010), 939--1008.

\bibitem{LOST11}  M.~Liebeck, E.~A.~O'Brien, A.~Shalev, Pham Huu Tiep,
 {\it Products of squares in finite simple groups}.
---  Proc. Amer. Math. Soc.  (2011).

\bibitem{LP}  M.~Liebeck, L. Pyber,
 {\it Finite linear groups and bounded generation}.
---  Duke Math. J.  {\bf107} (2001), 159--171.

\bibitem{liehl}  B.~Liehl,
{\it Beschr\"ankte Wortl\"ange in\/ $\SL_2$}.
 ---  Math. Z. {\bf186} (1984), 509--524.

\bibitem{matsumoto} H.~Matsumoto,
{\it Sur les sous-groupes arithm\'etiques des groupes semi-simples d\'eploy\'es}.
---  Ann.\ Sci.\ \'Ecole Norm.~Sup. (4) {\bf2} (1969), 1--62.

\bibitem{moree99}  P.~Moree,
 {\it On primes in arithmetic progression having a prescribed primitive root}.
---  J. Number Theory {\bf78} (1999), 85--98.

\bibitem{moree08}  P.~Moree,
 {\it On primes in arithmetic progression having a prescribed
primitive root. {\rm II}}.
---  Funct. Approx. Comment. Math. {\bf39} (2008), 133--144.

\bibitem{morris}  D.~W.~Morris,
{\it Bounded generation of\/ $\SL(n,A)$ {\rm(}after D.~Carter,
G.~Keller, and E.~Paige{\rm)}}
---  New York J. Math. {\bf13} (2007), 383--421.

\bibitem{NDS}  K.~R.~Nagarajan, M.~P.~Devaasahayam, T.~Soundararajan,
 {\it Products of three triangular matrices over commutative rings}.
---  Linear Algebra Applic. {\bf348} (2002), 1--6.

\bibitem{nikolov} N.~Nikolov,
 {\it A product decomposition for the classical quasisimple groups}.
---  J. Group Theory {\bf10} (2007), 43--53.

\bibitem{NP}  N.~Nikolov, L.~Pyber,
 {\it Product decomposition of quasirandom groups and
a Jordan type theorem}.
  arXiv:math/0703343 (2007).

\bibitem{paeth}  A.~Paeth,
 {\it A fast algorithm for general raster rotation}.
--- In: Graphics Gems, Acad. Press (1990), pp.~179--195.

\bibitem{RR}  A.~S.~Rapinchuk, I.~A.~Rapinchuk,
 {\it Centrality of the congruence kernel for elementary
sub\-groups of Chevalley groups of rank $>1$ over
Noetherian rings}. (2010), pp.~1--12
  arXiv:1007.2261v1 [math.GR].

\bibitem{shalev07} A.~Shalev,
 {\it Commutators, words, conjugacy classes, and character methods}.
---  Turk. J. Math. {\bf31} (2007), 131--148.

\bibitem{shalev09} A.~Shalev,
 {\it Word maps, conjugacy classes, and a noncommutative
Waring-type theorem}.
---  Ann. Math. {\bf170}, No.~3 (2009), 1383--1416.

\bibitem{sharpe} R.~W.~Sharpe
 {\it On the structure of the Steinberg group $\St(\Lambda)$}.
---  J.~Algebra {\bf68} (1981), 453--467.

\bibitem{SH} She Yiyuan, Hao Pengwei,
 {\it On the necessity and sufficiency of PLUS factorizations}.
---  Linear Algebra Applic. {\bf400} (2005), 193--202.

\bibitem{SiV}  S.~Sinchuk, N.~Vavilov,
 {\it Parabolic factorisations of exceptional Chevalley groups}
(to appear).

\bibitem{SS}  A.~Sivatski, A.~Stepanov,
{\it On the word length of commutators in\/ $\GL_n(R)$}.
 ---  $\K$-theory {\bf17} (1999), 295--302.

\bibitem{stein}  M.~R.~Stein,
 {\it Surjective stability in dimension\/ $0$ for\/ $\K_{2}$ and related
functors}.
 ---  Trans. Amer. Math. Soc. {\bf178} (1973), 176--191.

\bibitem{SV}  A.~Stepanov, N.~Vavilov,
{\it On the length of commutators in Chevalley groups}.
---  Israel Math. J. (2011), 1--20.

\bibitem{strang}  G.~Strang,
 {\it Every unit matrix is a $LULU$}.
---  Linear Algebra Applic. {\bf265} (1997), 165--172.

\bibitem{tavgen92}  O.~I.~Tavgen,
 {\it Bounded generation of normal and twisted Chevalley
groups over the rings of\/ $S$-integers}.
 ---  Contemp. Math. {\bf131} (1992), No.~1, 409--421.

\bibitem{toffoli}  T.~Toffoli,
 {\it Almost every unit matrix is a $ULU$}.
---  Linear Algebra Applic. {\bf259} (1997), 31--38.

\bibitem{TQ} T.~Toffoli, J.~Quick,
 {\it Three dimensional rotations by three shears}.
---  Graphical Models \& Image Processing {\bf59} (1997), 89--96.

\bibitem{vaser}  L.~N.~Vaserstein,
{\it Bass's first stable range condition}.
---  J. Pure Appl. Algebra {\bf34} (1984), Nos.~2--3, 319--330.

\bibitem{VW1}  L.~N.~Vaserstein, E.~Wheland,
{\it Factorization of invertible matrices over rings of
stable rank one}.
--- J. Austral. Math. Soc., Ser.~A, {\bf 48} (1990), 455--460.

\bibitem{VW2}  L.~N.~Vaserstein, E.~Wheland,
 {\it Commutators and companion matrices over rings of stable rank $1$}.
---  Linear Algebra Appl. {\bf142} (1990), 263--277.

\bibitem{vavilov}  N.~Vavilov,
{\it Structure of Chevalley groups over commutative rings}.
--- In:  Proc. Conf. Non-as\-so\-ci\-a\-ti\-ve algebras and related topics
(Hiroshima -- 1990) , World Sci. Publ., London
et al. (1991), pp.~219--335.

\bibitem{VP}  N.~Vavilov, E.~Plotkin,
{\it Chevalley groups over commutative rings. {\rm I.} Elementary
calcula\-tions}.
 ---  Acta Applicandae Math. {\bf45} (1996), 73--115.

\bibitem{you} You Hong, Commutators and unipotents in
symplectic groups. --- Acta Math. Sinica, New Ser., {\bf 10}
(1994), 173--179.

\bibitem{ZY} Zheng Bodong, You Hong, Products of commutators
of transvections over local rings. --- Linear Algebra Applications
{\bf 357} (2002), 45--57.

\end{thebibliography}
\end{document}